\documentclass{article}

\usepackage[margin=1in]{geometry}

\usepackage{graphicx} 

\usepackage{amsmath,amsfonts,amsthm,amssymb,mathrsfs,dsfont,bbm} 
\usepackage{mathtools}
\usepackage{bbm}
\usepackage{hyperref}
\usepackage{subcaption}
\usepackage[dvipsnames]{xcolor}
\usepackage{algorithm}
\usepackage{algpseudocode}
\usepackage{natbib}
\usepackage{dsfont}
\usepackage{enumitem}
\usepackage{comment}

\newcommand{\cX}{\mathcal{X}}

\newcommand{\cT}{\mathcal{T}}
\newcommand{\bP}{\mathbb{P}}
\newcommand{\bE}{\mathbb{E}}
\newcommand{\Z}{\mathbb{Z}}

\newcommand{\1}{\mathbbm{1}}

\newtheorem{theorem}{Theorem}[section]

\newtheorem{lemma}[theorem]{Lemma}
\newtheorem{proposition}[theorem]{Proposition}

\theoremstyle{definition}
\newtheorem*{remark}{Remark}

\renewcommand{\P}{\bP}
\newcommand{\E}{\bE}      

\DeclarePairedDelimiter\floor{\lfloor}{\rfloor}

\title{Shotgun Assembly of the Non-Uniform Labeled Lattice Model}
\author{Julia Gaudio, Charlie K. Guan, Andrew Jin}
\date{\today}

\begin{document}

\maketitle

\begin{abstract}
    In this paper, we study the graph shotgun assembly problem in the labeled lattice model, as introduced by Mossel and Ross (2019). We start with a $d$-dimensional lattice of length $n$, where the color of each vertex is sampled i.i.d. from a distribution $\pi$. Then, we ask whether it is possible to reconstruct the coloring of the lattice given the coloring of each subgrid of length $k$. While previous work by Ding and Liu (2023) and Demidovich, Panichkin, and Zhukovskii (2023) identified the reconstruction threshold when $\pi$ is uniform, we examine the non-uniform case. This paper establishes a sufficient condition for the impossibility of reconstruction for $d \geq 2$ using information-theoretic arguments.  Our negative result matches the known threshold in the uniform case. 
\end{abstract}

\section{Introduction}
\label{sec:intro}
The graph shotgun assembly problem aims to infer a graph when given collections of local neighborhoods of the graph. One motivating application is the reconstruction of neural connections in neuroscience: given sub-networks of neurons collected from experiments, the goal is to reconstruct the entire neural network of the brain \cite{soudry2015}. The term \textit{shotgun assembly} is inspired by the problem of DNA shotgun assembly, in which an entire DNA sequence is reconstructed from fragments of DNA subsequences that overlap each other. In this setting, we have a long 1-dimensional lattice of nodes, each labeled with one of four values (adenine, cytosine, guanine, or thymine), and we would like to recover the label of each node based on observations of smaller, overlapping fragments.

The labeled lattice model, proposed by Mossel and Ross \cite{MosselRoss2019}, generalizes DNA shotgun assembly for any arbitrary dimension and number of labels. Let $\pi$ be the probability mass function of a discrete distribution on $r$ values. Then, the labeled lattice model is defined as the $d$-dimensional lattice of length $n$, denoted $\Z_n^d$, with $\text{i.i.d.}$ vertex labels (also called colors) from a fixed $\pi$. For the rest of this paper, we will use ``lattice'' to refer to the labeled lattice model. 

The shotgun assembly problem on the labeled lattice model asks when it is possible to recover the coloring of the whole lattice, given that we can observe the coloring of small subgrids (but not their location within the lattice). To formalize this, we define a \textit{$k$-card} of the lattice as a subgrid of length $k$ where the coloring and orientation are observed, but the location is not. Then, we define the \textit{$k$-deck} of the lattice as the collection of all $k$-cards. Note that each $k$-card has $k^d$ vertices, and that the $k$-deck contains $(n-k+1)^d$ total cards. Then, we call a lattice \textit{identifiable} from its $k$-deck if no other lattice shares the same $k$-deck, and \textit{non-identifiable} otherwise. The shotgun assembly problem aims to 1) determine the values of $k$ for which the lattice is identifiable and non-identifiable, and 2) if the lattice is identifiable, find an efficient algorithm to recover the lattice coloring from the $k$-deck.

Prior work has largely focused on the case where $\pi$ is a uniform discrete distribution. The threshold for $d=1$ was found by Dyer, Frieze, and Suen \cite{Dyer_Frieze_Suen}, who showed that as $n \to \infty$, the lattice is identifiable with a constant nontrivial probability if $k = \frac{2 \log n - \log(c)}{\log r}$ where $c$ is a constant. Consequently, the lattice coloring is identifiable with high probability when $k \geq \frac{2(1+\epsilon) \log n}{\log r}$ while it is non-identifiable with high probability when $k \leq \frac{2(1-\epsilon) \log n}{\log r}$.

The threshold for $d \geq 2$ was first found by Ding and Liu \cite{ding2023lattice}, who showed that if $k^d \geq \frac{(1+\epsilon) \log n^d}{\log r}$, then the lattice is identifiable with high probability, while if $k^d \leq \frac{(1-\epsilon) \log n^d}{\log r}$, the lattice is non-identifiable with high probability. We note that Ding and Liu also developed an efficient algorithm to recover the lattice coloring above the threshold. Later, Narayanan and Yap \cite{narayanan2025reconstructing} improved the threshold to two-point concentration in the case where $d=2$ and $r = 2$ using an alternative algorithm, which was later generalized to all $d\geq2$ and $r \geq 2$ by Demidovich, Panichkin, and Zhukovskii \cite{demidovich2023}. The latter showed that if $k \geq (\frac{\log n^d }{\log r})^{1/d} + \frac{1}{d} + \epsilon$, then the lattice is identifiable with high probability, while if $k \leq (\frac{\log n^d}{\log r})^{1/d} - \epsilon$, then the lattice is non-identifiable with high probability. Thus, the threshold for identifiability is well-understood in the uniform case.

Moving to the case where $\pi$ is non-uniform, the threshold for $d=1$ was also identified by Dyer, Frieze, and Suen \cite{Dyer_Frieze_Suen}; the authors showed that as $n \to \infty$, the lattice can be recovered with a constant nontrivial probability when $k = \frac{2 \log n - \log c}{\log(1/q)}$ where $q := \sum_{i=1}^r \pi_i^2$ and $c$ is a constant. Thus, the lattice coloring is identifiable with high probability when $k \geq \frac{2(1+\epsilon) \log n}{\log(1/q)}$, and non-identifiable with high probability when $k \leq \frac{2(1-\epsilon) \log n}{\log(1/q)}$.

However, much less is known for $d \geq 2$ since the analysis becomes substantially more difficult. Mossel and Ross \cite{MosselRoss2019} identified a sufficient condition for identifiability in which all $(k-1)$-cards are unique, which is $k^d \geq \frac{2(1+\epsilon) \log n^d}{\log(1/q)}$. Under this condition, we can recover the lattice by starting with an arbitrary initial $k$-card, and then looking for other cards with a shared subgrid of length $k-1$; since all subgrids with length $k-1$ are unique in the lattice, cards with shared subgrids must necessarily be adjacent. We note that their sufficient condition does not match the threshold in the uniform case. To the best of our knowledge, this was the only known result in the non-uniform case.

Our main contribution is establishing a sufficient condition for non-identifiability in the non-uniform case, which matches the known threshold in the uniform case for $d \geq 2$. We conjecture that our sufficient condition for non-identifiability is the true threshold in the non-uniform case. In addition, we also provide an alternative proof of non-identifiability for the $d=1$ case, which could help with understanding the lattice's properties better.

\subsection{Related Work}
There is also a large body of work on related shotgun assembly problems that were first introduced by Mossel and Ross \cite{MosselRoss2019}, and we discuss some of them here. 

\paragraph{Dense Erd{\H{o}}s--R{\'e}nyi}
One of the most well-studied such problems is shotgun assembly of unlabeled Erd{\H{o}}s--R{\'e}nyi graphs, which asks whether we can reconstruct a graph (up to isomorphism) from the local $r$-neighborhood of each vertex. In the dense $G(n, p_n)$ setting, the focus has largely been on finding the edge probabilities $p_n$ such that reconstruction is possible for a fixed $r$. There are two phase transitions---reconstruction is possible when the probability is very small or very large, and impossible in between. Intuitively, when $p_n$ is very small, all components fit within a single $r$-neighborhood, while if $p_n$ is very large, all $r$-neighborhoods are unique \cite{Johnston2025}. 

For $r=1$, work by Gaudio and Mossel \cite{gaudioMossel2022}, Huang and Tikhomirov \cite{HuangTikhomirov2025}, and Johnston, Kronenberg, Roberts, and Scott \cite{Johnston2025} have identified the phase transitions at $p_n = \Theta(n^{-3/2})$ and $p_n = \widetilde{\Theta}(n^{-1/2})$, where $\widetilde{\Theta}$ denotes equality up to polylog factors.

For $r=2$, Johnston et al. \cite{Johnston2025} identified one phase transition at $p_n = \Theta(n^{-5/4})$. The other phase transition is still unknown; Gaudio and Mossel \cite{gaudioMossel2022} showed that if $p_n \leq n^{-3/4-\epsilon}$ then reconstruction is impossible, while Johnston et al. \cite{Johnston2025} showed that if $p_n \geq n^{-2/3+\epsilon}$ then reconstruction is possible.

For $r=3$, the phase transitions are known from Johnston et al. \cite{Johnston2025}, expanding on the initial work of Mossel and Ross \cite{MosselRoss2019}. One transition occurs at $p_n = \Theta(n^{-\frac{2r+1}{2r}})$, while the other transition occurs at $p_n = \Theta(\frac{\log^2 n}{n (\log \log n)^3})$.

Finally, for $r \geq 4$, Johnston et al. \cite{Johnston2025} identified the phase transitions at $p_n = \Theta(n^{-\frac{2r+1}{2r}})$ and $p_n = \Theta(\frac{\log n}{n})$. A more precise result is known for $r = 4$: Gaudio, R{\'a}cz, and Sridhar \cite{gaudio2025average} showed that reconstruction is possible when $p_n \geq (1+\epsilon)\frac{\log n}{n}$.

\paragraph{Sparse Erd{\H{o}}s--R{\'e}nyi} The sparse setting considers whether reconstruction is possible in $G(n, \lambda/n)$, which depends on the value of $r$. Initial work by Mossel and Ross \cite{MosselRoss2019} showed that the threshold occurs at $r = \Theta(\log n)$. Later, Ding, Jiang, and Ma \cite{DingJiangYiMa2023} identified a sharp threshold at $r = \frac{\log n}{\log(\lambda^2 \gamma_\lambda)^{-1}}$, where $\gamma_\lambda$ is the probability that two independent Galton-Watson trees with parameter $\lambda$ are isomorphic. In particular, they showed that if $r \leq \frac{(1-\epsilon) \log n}{\log(\lambda^2 \gamma_\lambda)^{-1}}$ then reconstruction is impossible, while if $r \geq \frac{(1+\epsilon) \log n}{\log(\lambda^2 \gamma_\lambda)^{-1}}$ then reconstruction is possible via a polynomial time algorithm.

\paragraph{Random Jigsaw Puzzles} 
Another related problem introduced by Mossel and Ross \cite{MosselRoss2019} is reconstructing random ``jigsaw puzzles''. In this problem, we consider a grid of $n \times n$ squares, where each border between adjacent squares has one of $q$ configurations drawn uniformly at random. Then, we want to know the values of $q$ for which we can reassemble the squares back together. 

Bordenave, Feige, and Mossel \cite{bordenave2020shotgun} and Nenadov, Pfister, and Steger \cite{nenadov2016unique} initially showed that if $q \geq n^{1+\epsilon}$, then the puzzle has a unique solution, while if $q = o(n)$, then the puzzle does not have a unique solution. Follow-up work by Balister, Bollob{\'a}s, and Narayanan \cite{balister2017reconstructing} showed that there exist constants $c, C$ such that $q \geq Cn$ implies that the puzzle is unique, while $q \leq cn$ implies that the puzzle is not unique. Then, Martinsson \cite{martinsson2019linear} identified specific constants, showing that if $q \geq (2+\epsilon)n$ then the puzzle is unique, while if $q \leq \frac{2}{\sqrt{e}}n$ then the puzzle is not unique. The true threshold remains an open problem.

\paragraph{}
The remainder of the paper is organized as follows. In Section \ref{sec:main_results}, we present our main result, a sufficient condition for non-identifiability when $d \geq 2$. We subsequently prove our main result in Section \ref{sec:2d_proof}. In Section \ref{sec:1d_proof}, we provide an alternative proof of non-identifiability for $d=1$. In Section \ref{sec:threshold_discussion}, we discuss potential candidates for the true threshold for $d \geq 2$. Finally, we conclude with some future research directions in Section \ref{sec:conclusion}.

\section{Main Result}
\label{sec:main_results}
In this section, we present a sufficient condition for the lattice to be non-identifiable, which matches the known threshold in the uniform case \cite{ding2023lattice}.

\begin{theorem}
\label{thm:main_theorem_2d_negative}
    Let $d \geq 2$ and fix $\epsilon > 0$. Then, the lattice is non-identifiable with high probability if
    \begin{equation}
    \label{eq:2d_negative_condition}
        k^d \leq \frac{(1-\epsilon) \log n^d}{H(\pi)},
    \end{equation}
    where $H(\pi) := - \sum_{i=1}^r \pi_i \log (\pi_i)$ denotes the entropy of $\pi$.
\end{theorem}
\begin{remark}
    Observe that if $\pi$ is uniform, then $H(\pi) = \sum_{i=1}^r \frac{1}{r} \log r = \log r$. Thus, our condition for non-identifiability becomes $k^d \leq \frac{(1-\epsilon) \log n^d}{\log r}$, which is indeed the threshold in the uniform case.
\end{remark}

We use techniques from information theory to prove Theorem \ref{thm:main_theorem_2d_negative}. In particular, we use Fano's inequality to lower bound the probability that the MAP estimator fails to recover the lattice from the deck, which shows that the error probability is $1 - o(1)$. Intuitively, under Condition \eqref{eq:2d_negative_condition}, the entropy of the lattice exceeds the entropy of the $k$-deck, so there is not enough information in the $k$-deck to recover the lattice.

\section{Proof of Non-identifiability for $d \geq 2$}
\label{sec:2d_proof}
In this section, we show that the lattice is non-identifiable with high probability when $k^d \leq \frac{(1-\epsilon)\log n^d}{H(\pi)}$.

\subsection{Preliminaries}
We first state several lemmas that will be used in the proof.

\begin{lemma}[Fano's Inequality (Theorem 2.11.1 in \cite{CoverThomasM2006Eoit})]
    Let $X$, $Y$ be correlated random variables. Suppose that we observe $Y$, and we want to estimate $X$ from the observed value of $Y$, which we denote $\widehat{X}$. Define the error probability
    \begin{equation*}
        P_e = \P(\widehat{X} \neq X).
    \end{equation*}
    Then, for any estimator $\widehat{X} = g(Y)$
    \begin{equation}
        P_e \geq \frac{H(X|Y)-1}{\log |\cX|},
    \end{equation}
    where $H(X | Y)$ denotes the conditional entropy of $X$ conditioned on $Y$, and $\cX$ is the support of $X$.
\end{lemma}

\begin{lemma}[Fannes-Audenaert Inequality \cite{audenaert2007sharp}]
\label{lem:Fannes-Audenaert}
    Let $P$ and $Q$ be discrete random variables on the same support of $N$ elements with total variation distance $d_{TV}(P,Q) < \epsilon$. Then 
    \begin{equation*}
        |H(P)-H(Q)| \leq H(\epsilon) + \epsilon \log N.
    \end{equation*}
    where $H(\epsilon) := - \epsilon \log \epsilon - (1-\epsilon) \log (1-\epsilon)$ is the binary entropy of $\epsilon$.
\end{lemma}

\begin{lemma}[Hoeffding's Inequality]
    Let $X \sim \mathrm{Bin}(n, p)$. Then, 
    \begin{equation*}
        \P(|X - \E[X]| \geq t) \leq 2 \exp\bigg(-\frac{2t^2}{n}\bigg).
    \end{equation*}
\end{lemma}

\begin{lemma}[Binomial Chernoff Bound]
    Let $X \sim \mathrm{Bin}(n, p)$ and $\mu := \E[X]$. Then,
    \begin{equation*}
        \P(X < (1 - \delta)\mu) \leq \exp\bigg(-\frac{\mu\delta^2}{2}\bigg).
    \end{equation*}
\end{lemma}

Let $G$ be a lattice drawn from the labeled lattice model, and let $D$ be the corresponding $k$-deck. To show that the lattice is non-identifiable, we can show that any estimator $\widehat{G}$ fails to recover $G$ with high probability, meaning that $\P(\widehat{G} \neq G) = 1 - o(1)$. Rather than analyzing an estimator $\widehat{G}$ over all possible lattices, we instead consider a simpler set of lattices whose probability goes to 1 as $n \to \infty$, and show that any estimator fails to recover such lattices with high probability.

Fix a sequence $\delta_n = n^{-d/3}$, which implies that $\delta_n^2 n^d \to \infty$ and that $\delta_n \to 0$ as $n \to \infty$. We define $T_{\delta_n}$ as the event that the $n$-lattice has a ``typical'' coloring $\sigma$, in the sense that for all $i\in[r]$, 
\begin{equation*}
    \big|\{v: \sigma(v)=i\}\big| \in \Big[(1-\delta_n)\pi_i n^d, (1+\delta_n)\pi_i n^d\Big].
\end{equation*}
That is, the number of vertices of each color is close to their expectation. Thus, we call a lattice \textit{$\delta_n$-typical} if $T_{\delta_n}$ holds. It turns out that the lattice coloring is typical with high probability.

\begin{lemma}
\label{lem:typical_lattice_whp}
    The event $T_{\delta_n}$ holds with high probability.
\end{lemma}
\begin{proof}
    Let $\pi_\text{min} := \min_i \pi_i$. Observe that the number of vertices of color $i$ follows a $\text{Bin}(n^d, \pi_i)$ distribution. By Hoeffding's inequality and the union bound, 
    \begin{equation}
    \label{eq:typical_lattice_prob}
        \bP(T_{\delta_n}) \geq 1-  \sum_{i\in[r]}  2\exp(-2\delta_n^2 \pi_i^2 n^d) \geq 1-2r \exp(-2\delta_n^2 n^d)= 1- o(1),
    \end{equation}
    which shows that $T_{\delta_n}$ holds with high probability.
\end{proof}

Since the lattice coloring is typical with high probability, it suffices to show that typical lattices are not identifiable with high probability. To this end, we define $\widetilde{G}$ as a lattice coloring with the same distribution as $G$ conditioned on $T_{\delta_n}$, and $\widetilde{D}$ as a $k$-deck with the same distribution as $D$ conditioned on $T_{\delta_n}$. Let $G'$ be the MAP estimator for estimating $\widetilde{G}$ from $\widetilde{D}$. Then, showing that typical lattices are non-identifiable with high probability is equivalent to showing that $\P(G' \neq \widetilde{G}) = 1- o(1)$.

In addition, we define $\cT_{\delta_n}$ as the set of lattice colorings that satisfy $T_{\delta_n}$, and denote $M_n=|\cT_{\delta_n}|$, which is the number of $\delta_n$-typical $n$-lattices.

\begin{proposition}
\label{prop:2d_negative_fano}
     Suppose that the following three properties hold.
    \begin{enumerate}
        \item $H(\widetilde{G}) = n^d H(\pi) + o(1)$,
        \item $\log M_n \leq n^d H(\pi) + o(n^d)$, and
        \item $H(\widetilde{D}) = o(n^d)$.
    \end{enumerate}
    Then, $P(G' \neq \widetilde{G}) = 1 - o(1)$, which shows that the lattice is non-identifiable with high probability.
\end{proposition}
\begin{proof}
    From Fano's inequality, we have that 
    \begin{equation*}
        P(G' \neq \widetilde{G}) \geq \frac{H(\widetilde{G} \mid \widetilde{D}) - 1}{\log M_n} = \frac{H(\widetilde{G}) - H(\widetilde{D}) - 1}{\log M_n},
    \end{equation*}
    where $H(\widetilde{G} \mid \widetilde{D}) = H(\widetilde{G}) - H(\widetilde{D})$ because $\widetilde{D}$ is completely determined by $\widetilde{G}$. Then, using Properties 1-3 above, we obtain
    \begin{equation*}
        P(G' \neq \widetilde{G}) \geq \frac{n^dH(\pi) + o(1) - o(n^d) - 1}{n^d H(\pi) + o(n^d)} = 1 - o(1).
    \end{equation*}
\end{proof}

In the sequel, we prove that the three properties in Proposition \ref{prop:2d_negative_fano} hold when $k^d \leq \frac{(1-\epsilon) \log n^d}{H(\pi)}$, which allows us to conclude Theorem \ref{thm:main_theorem_2d_negative}.

\subsection{Proof of Property 1}
\begin{lemma}[Property 1]
\label{lem:property_1}
    We have that $H(\widetilde{G}) = n^d H(\pi) + o(1)$.
\end{lemma}
\begin{proof}
    Observe that we can represent the lattice as $n^d$ i.i.d. random variables which encode the color of each vertex in the lattice. Since the joint entropy of i.i.d. random variables equals the sum of their individual entropies, we obtain that $H(G) = n^d H(\pi)$. Hence, it is sufficient to show that $H(\widetilde{G}) = H(G) + o(1)$.

    To accomplish this, we first compute the total variation distance between $G$ and $\widetilde{G}$, denoted $d_{TV}(G, \widetilde{G})$. By the coupling inequality (Lemma 4.1.11 in \cite{Roch2024}), we have that
    \begin{equation*}
        d_{TV}(G, \widetilde{G}) \leq \P(G \neq \widetilde{G})
    \end{equation*}
    for all couplings of $G$ and $\widetilde{G}$. Thus, we consider the following coupling: If $G \in T_\delta$ then we let $\widetilde{G} = G$. If $G \notin T_\delta$, then we sample $\widetilde{G}$ separately so that the marginal distribution of $\widetilde{G}$ is preserved. In this case, we have 
    \begin{equation}
    \label{eq:G_total_variation_bound}
         d_{TV}(G, \widetilde{G}) \leq \P(G \neq \widetilde{G}) = 1 - \P(T_{\delta_n}) \leq 2r \exp(-2\delta_n^2 n^d)
    \end{equation}
    from \eqref{eq:typical_lattice_prob}. Now, we apply the Fannes-Audenaert inequality, which yields
    \begin{equation*}
        |H(G) - H(\widetilde{G})| \leq H\big(2r \exp(-2\delta_n^2 n^d)\big) + 2r \exp(-2\delta_n^2 n^d) \log r^{n^d}.
    \end{equation*}
    To analyze the first term, recall that $\delta_n^2 n^d \to \infty$ as $n \to \infty$. Therefore, we have that $2r \exp(-2\delta_n^2 n^d) \to 0$, which implies that $H(2r\exp(-2\delta_n^2 n^d)) \to 0$ as well. For the second term, we use the fact that $\delta = n^{-d/3}$ to obtain
    \begin{equation*}
        2r \exp(-2\delta_n^2n^d) \log r^{n^d} = 2r  \exp(-2n^{d/3}) n^d \log r\to 0.
    \end{equation*}
    Therefore, $|H(G) - H(\widetilde{G})| = o(1)$, which shows that $H(\widetilde{G}) = H(G) + o(1) = n^dH(\pi) + o(1)$.
    
\end{proof}

\subsection{Proof of Property 2}
\begin{lemma}[Property 2]
    We have that $\log M_n \leq n^d H(\pi) + o(n^d)$.
\end{lemma}
\begin{proof}
    Observe that the probability of any given $\delta_n$-typical lattice is at least 
    \begin{equation*}
        \prod_{i=1}^r \pi_i^{(1+\delta_n) \pi_i n^d}
    \end{equation*}
    since there are at most $(1+\delta_n) \pi_i n^d$ vertices with color $i$. Therefore, the number of typical lattices is at most
    \begin{equation*}
        M_n \leq \bigg(\prod_{i=1}^r \pi_i^{(1+\delta_n) \pi_i n^d}\bigg)^{-1}
    \end{equation*}
    which implies that 
    \begin{equation*}
        \log M_n \leq - \sum_{i=1}^r (1+\delta_n)\pi_i n^d \log (\pi_i) = (1+\delta_n)n^d \sum_{i=1}^r -\pi_i \log \pi_i = (1+\delta_n)n^d H(\pi).
    \end{equation*}
    Since $\delta_n \to 0$, we have that $\log M_n \leq (1 + o(1))n^d H(\pi) = n^dH(\pi) + o(n^d)$, which completes the proof.
\end{proof}

\subsection{Proof of Property 3}
We prove Property 3 using a similar strategy as Property 1. First, we show that $H(D) = o(n^d)$ by finding an explicit encoding for $D$ and showing that the expected encoding length is $o(n^d)$; since $H(D)$ is less than the expected encoding length of any encoding \cite{CoverThomasM2006Eoit}, it follows that $H(D) = o(n^d)$ as well. Then, we use the coupling inequality to compute the total variation distance between $D$ and $\widetilde{D}$, and apply the Fannes-Audenaert inequality to show that $|H(D) - H(\widetilde{D})| = o(1)$.

\begin{lemma}
\label{lem:entropy_of_deck}
    If $k^d \leq \frac{(1-\epsilon) \log n^d}{H(\pi)}$, then $H(D) = o(n^d)$.
\end{lemma}
\begin{proof}
We construct an explicit encoding of the deck and show that the expected encoding length is $o(n^d)$. Consider the following encoding: For every type of card that appears in the deck, we will encode the type of card using a unique identifier of length $\log(r^{k^d})$, then encode the number of times that type of card appears using a string of length $\log(n^d)$. 

Observe that there are at most $r^{k^d}$ types of cards because there are $k^d$ vertices in each card and $r$ choices for the color of each vertex. Thus, we can give each type of card a unique identifier using $\log(r^{k^d})$ bits. Meanwhile, there are less than $n^d$ cards in the deck (which is also the maximum number of times any given type of card can appear), so we can encode the count of each card type using $\log(n^d)$ bits.

Now, we compute the expected encoding length of the deck from the point of view of a fixed, unrevealed card. Its contribution to the encoding length is given by
\[\frac{1}{Z} \left(\log(r^{k^d}) + \log(n^d) \right),\]
where $Z$ is the number of times that card appears in the deck---this is an overcounting factor, since all cards of the same type have a collective encoding length of $\log(r^{k^d}) + \log(n^d)$. Thus, the expected encoding length of the deck is
\begin{equation*}
    \mathbb{E}\left[\frac{1}{Z}\right] \cdot n^d \left(\log(r^{k^d}) + \log(n^d) \right).
\end{equation*} 
It remains to bound $\mathbb{E}[1/Z]$. First, observe that for any $t > 0$,
    \begin{align}
    \label{eq:expectation_1/Z}
        \mathbb{E}\left[\frac{1}{Z}\right] &= \mathbb{E}\left[\frac{1}{Z} ~\bigg|~ Z \geq t \right] \mathbb{P}(Z \geq t) + \mathbb{E}\left[\frac{1}{Z} ~\bigg|~ Z < t \right] \mathbb{P}(Z < t) \notag \\
        &\leq \frac{1}{t} \mathbb{P}(Z \geq t) + \mathbb{P}(Z < t) \notag \\
        &\leq \frac{1}{t} + \mathbb{P}(Z < t).
    \end{align}
    To upper-bound $\E[1/Z]$, we let $t = n^{\epsilon d / 2}$ and focus on upper-bounding $\P(Z < t)$. Rather than examine $Z$ directly, we define a new random variable $Z'$ as follows: Partition the lattice into non-overlapping regions of $k^d$ vertices, which we call the cards of the \textit{partitioned lattice}; there are a total of $(n/k)^d$ cards in the partitioned lattice. Then, consider a random card from the partitioned lattice, which we call the top card, and let $Z'$ be the number of times the top card appears in the partitioned lattice. We note that $Z' \preceq Z$ stochastically, which yields $\bP(Z < t) \leq \bP(Z' < t)$.
    
    We now upper-bound $\bP(Z' < t)$. Let $X = (X_1, \ldots, X_r)$, where $X_i$ is the number of vertices of color $i$ in the top card. For any fixed $\delta > 0$, define $E_\delta$ as the event
    \begin{equation}
    \label{eq:E_delta}
        E_\delta = \Big\{X_i \in\big[(1-\delta)\pi_i k^d, (1+\delta)\pi_i k^d\big] \: \forall i \in [r]\Big\},
    \end{equation}
    meaning that the number of vertices of each color is close to its expectation. We call the top card \textit{$\delta$-typical} if $E_\delta$ holds. 
    Note that the top card is $\delta$-typical with high probability, which can be shown using a similar argument to Lemma \ref{lem:typical_lattice_whp}.

    To compute $\P(Z' < t)$, we fix $\delta = \epsilon/4$ and use the law of total probability, conditioning on $X$. Furthermore, we divide the sum into two parts depending on whether the card is $\delta$-typical, which yields  
    \begin{align}
    \label{eq:probability_Z'}
    \P(Z' < t) &= \sum_{x \in E_\delta} \P(Z' < t \mid X = x) \P(X=x) + \P\Big(\{Z' < t\} \cap \{X \notin E_\delta\} \Big) \notag \\
    &= \sum_{x \in E_\delta} \P(Z' < t \mid X = x) \P(X=x) + o(1).
    \end{align}
    Conditioned on $X = x$, the probability that any given card equals the top card is $\prod_{i=1}^r \pi_i^{x_i}$. Therefore, the number of cards matching the top card is distributed according to $(Z' \mid X = x) \sim \text{Bin}((n/k)^d, \prod_{i=1}^r \pi_i^{x_i})$. Now, define $Y \sim \text{Bin}((n/k)^d, \prod_{i=1}^r \pi_i^{(1 + \delta)\pi_i k^d})$, and since $x \in E_\delta$, we have that $x_i \leq (1+\delta)\pi_ik^d$ for all $i$. Hence, $Y \preceq (Z' \mid X = x)$ stochastically for any $x \in E_\delta$, which implies that $\P(Z' < t \mid X = x) \leq \P(Y < t)$ for any $x \in E_\delta$. In light of \eqref{eq:probability_Z'}, we obtain that $\P(Z' < t) \leq \P(Y < t)$, which also implies that 
    \begin{equation}
    \label{eq:compare_prob_Y_Z}
        \P(Z < t) \leq \P(Y < t).
    \end{equation}
    
    Now, we compute an upper bound for $\P(Y < t)$ using the Chernoff bound. We first compute $\E[Y]$, which yields \begin{align*}
        \E[Y] = \frac{n^d}{k^d} \prod_{i=1}^r \pi_i^{(1+\delta)\pi_ik^d} &= \frac{n^d}{k^d} \exp\bigg((1+\delta)k^d\sum_{i=1}^r \pi_i \log \pi_i \bigg) \\
        &= \frac{n^d}{k^d} e^{-(1+\delta)k^d H(\pi)}.
    \end{align*}
    Under the condition $k^d \leq \frac{(1-\epsilon)d \log n}{H(\pi)}$, we have that 
    \begin{align*}
        \frac{n^d}{k^d} \geq \frac{n^d H(\pi)}{(1-\epsilon)d \log n} = \frac{n^d}{C \log n},
    \end{align*}
    where $C := \frac{(1-\epsilon)d}{H(\pi)}$ is a constant. Under this condition, we also have that
    \begin{align*}
        e^{-(1+\delta)H(\pi)k^d} \geq e^{-(1+\delta)(1-\epsilon) d \log n} = n^{-(1+\delta)(1-\epsilon)d} > n^{-(1-3\epsilon/4)d},
    \end{align*}
    where the last inequality follows from $\delta = \epsilon/4$. Therefore, 
    \begin{align*}
        \E[Y] \geq \frac{n^d}{C \log n} n^{-(1-3\epsilon/4)d} = \frac{n^{3\epsilon d/4}}{C \log n} \gtrsim n^{ 2\epsilon d / 3}.
    \end{align*}
    Now, we use the Chernoff bound to upper-bound $\P(Y < t)$ where $t = n^{\epsilon d / 2}$, which yields
    \begin{align*}
    \label{eq:Chernoff_bound_Y}
        \P\Big(Y \leq n^{\epsilon d / 2}\Big) \leq \P\Big(Y \leq \frac{1}{2}n^{2\epsilon d/3}\Big) \leq \P\Big(Y \leq \frac{1}{2}\mu\Big) \leq \exp\big(-\frac{1}{8}\mu \big) \leq \exp\Big(-\frac{1}{8} n^{2\epsilon d / 3}\Big).
    \end{align*}
    From \eqref{eq:compare_prob_Y_Z}, we see that $\P(Z \leq n^{\epsilon d / 2}) \leq \exp(-\frac{1}{8} n^{2\epsilon d / 3})$ as well. Substituting this into \eqref{eq:expectation_1/Z}, we obtain that
    \begin{equation*}
        E\left[\frac{1}{Z}\right] \leq n^{-\epsilon d / 2} + \exp(-n^{2\epsilon d / 3} / 8) \leq n^{-\epsilon d / 3}.
    \end{equation*}
    Consequently, the expected encoding length of the deck is at most
    \begin{equation*}
        n^{d(1-\epsilon / 3)} \left(\log(r^{k^d}) + \log(n^d) \right) \lesssim  n^{d(1 -\epsilon/4)} = o(n^d),
    \end{equation*}
    which shows that $H(D) = o(n^d)$.
\end{proof}

Using Lemma \ref{lem:entropy_of_deck}, we can now prove Property 3. 

\begin{lemma}[Property 3]
    If $k^d \leq \frac{(1-\epsilon)\log n^d}{H(\pi)}$, then $H(\widetilde{D}) = o(n^d)$.
\end{lemma}
\begin{proof}
We first use the coupling inequality (Lemma 4.1.11 in \cite{Roch2024}) to bound the total variation distance between $D$ and $\widetilde{D}$, which yields
\begin{equation*}
    d_{TV}(D, \widetilde{D}) \leq \P(D \neq \widetilde{D}) 
\end{equation*}
for any coupling of $D, \widetilde{D}$. Consider the coupling between $G, \widetilde{G}$ described in Lemma \ref{lem:property_1}; since $D, \widetilde{D}$ are determined by $G, \widetilde{G}$ respectively, this also induces a coupling between $D, \widetilde{D}$. Furthermore, we note that if $G = \widetilde{G}$, then $D = \widetilde{D}$ as well. Therefore, we obtain 
\begin{equation*}
    d_{TV}(D, \widetilde{D}) \leq \P(D \neq \widetilde{D}) \leq  \P(G \neq \widetilde{G}) \leq 2r \exp(-2\delta_n^2n^d),
\end{equation*}
where the last inequality follows from \eqref{eq:G_total_variation_bound}. Then, applying the Fannes-Audenaert inequality and repeating the argument used in Lemma \ref{lem:property_1} shows that $|H(D) - H(\widetilde{{D}})| = o(1)$. Since $H(D) = o(n^d)$ from Lemma \ref{lem:entropy_of_deck}, we obtain that $H(\widetilde{D}) = o(n^d)$ as well.
\end{proof}

\section{Proof of Non-identifiability for $d=1$}
\label{sec:1d_proof}
In this section, we provide an alternative proof of non-identifiability in the 1-dimensional case. We first state the threshold established by Dyer, Frieze, and Suen \cite{Dyer_Frieze_Suen} and Mossel and Ross \cite{MosselRoss2019}.

\begin{theorem}
\label{thm:1d_threshold}
    Let $d = 1$ and fix $\epsilon > 0$. Then, the lattice is identifiable with high probability if
    \begin{equation}
    \label{eq:1d_positive_condition}
        k \geq \frac{2(1+\epsilon) \log n}{\log(1/q)},
    \end{equation}
    and non-identifiable with high probability if
    \begin{equation}
    \label{eq:1d_negative_condition}
        k \leq \frac{2(1-\epsilon) \log n}{\log(1/q)},
    \end{equation}
    where $q := \sum_{i=1}^r \pi_i^2$ is the probability that two vertices have the same color.
\end{theorem}

We note that this threshold corresponds to all $(k-1)$-cards being unique. Above this threshold, all $(k-1)$-cards are unique with high probability; this allows us to identify adjacent $k$-cards by looking for subgrids of length $k-1$ which appear in both cards, and then reconstruct the lattice by piecing together adjacent cards. Meanwhile, below that threshold, many $k$-cards are not unique, which allows us to produce multiple distinct lattices with the same $k$-deck through a swapping procedure.

To prove non-identifiability under Condition \eqref{eq:1d_negative_condition}, we extend a technique from Ding and Liu \cite{ding2023lattice} that showed non-identifiability in the uniform case. Let $I_1, \ldots, I_6$ be disjoint, consecutive intervals in the 1-dimensional lattice $\Z_n^1$ with $m := \floor{n/6}$ vertices each. For each $I_j$, let $\Gamma_j$ be a collection of $\ell := \floor{m/k}$ disjoint $k$-cards within $I_j$. We claim the following proposition holds. 

\begin{proposition}
\label{prop:1d_prop}
    Suppose that $k \leq \frac{2(1-\epsilon)\log n}{\log(1/q)}$. With high probability, there exist $B_j \in \Gamma_j$ for $j = 1, 3, 4, 6$ such that
    \begin{equation}
    \label{eq:non-identifiability_DingLiu_1d_condition1}
        \sigma|_{B_1} = \sigma|_{B_4} \text{ and } \sigma|_{B_3} = \sigma|_{B_6}
    \end{equation}
    and that
    \begin{equation}
    \label{eq:non-identifiability_DingLiu_1d_condition2}
        \sigma|_{J[1, m]} \neq \sigma|_{J'[1, m]},
    \end{equation}
    where $J$ is the interval strictly between $B_1$ and $B_3$ and $J'$ is the interval strictly between $B_4$ and $B_6$, and $J[1,m]$ is the first $m$ vertices in $J$ and $J'[1,m]$ is the first $m$ vertices in $J'$.
\end{proposition}

Proposition \ref{prop:1d_prop} is sufficient to show non-identifiability because $J$ and $J'$ can be swapped to produce different lattices while keeping the same deck; we refer the reader to Section 2.1 in \cite{ding2023lattice} for full details. Therefore, we obtain the non-identifiability result in Theorem \ref{thm:1d_threshold} as a corollary.

\begin{proof}[Proof of Proposition \ref{prop:1d_prop}]
We first show \eqref{eq:non-identifiability_DingLiu_1d_condition1} using the second moment method. Let $Z$ be the number of pairs $B \in \Gamma_1, B' \in \Gamma_4$ such that $\sigma|_{B} = \sigma|_{B'}$. Formally,
\begin{equation}
    Z = \sum_{B \in \Gamma_1, B' \in \Gamma_4} \1_{\sigma|_{B} = \sigma|_{B'}}.
\end{equation}
We will show that $\mathbb{P}(Z = 0) \to 0$ as $n \to \infty$ using Chebyshev's inequality. Observe that for any disjoint cards $B$ and $B'$,
\begin{equation*}
    \P(\sigma|_B = \sigma|_{B'}) = q^k    
\end{equation*}
because the probability that any two given vertices match is $q$. Hence,
\begin{equation*}
    \E[Z] = \ell^2 \bigg(\sum_{i=1}^r \pi_i^2\bigg)^k = \ell^2 q^k.
\end{equation*}
Now we compute the second moment. First, observe that
\begin{align}
\label{eq:non-identifiability_second_moment_initial}
    \E[Z^2] &= \E\bigg[\sum_{B \in \Gamma_1, B' \in \Gamma_4, C \in \Gamma_1, C' \in \Gamma_4} \1_{\sigma|_{B} = \sigma|_{B'}} \1_{\sigma|_{C} = \sigma|_{C'}} \bigg] \notag \\
    &= \ell^2 q^k + 2\ell^2(\ell-1) \bigg(\sum_{i=1}^r \pi_i^3\bigg)^k + \ell^2(\ell-1)^2 q^{2k} \notag \\
    &\leq \ell^2 q^k + 2\ell^3 \bigg(\sum_{i=1}^r \pi_i^3\bigg)^k + \ell^4 q^{2k}.
\end{align}
We will explain where each term in the second equality comes from. The first term comes from the case where $B = C$ and $B' = C'$. There are $\ell^2$ ways to choose such $B, B', C, C'$. To obtain $\sigma|_{B} = \sigma|_{B'}$ and $\sigma|_{C} = \sigma|_{C'}$, we need one subinterval of $k$ vertices to have the same coloring as another subinterval of $k$ vertices, which occurs with probability $q^k$. 

The second term comes from the case where either $B = C$ or $B' = C'$, but not both. Assume for now that $B = C$ and $B' \neq C'$. There are $\ell \cdot \ell(\ell-1) = \ell^2(\ell-1)$ ways to choose such $B, B', C, C'$ (the $\ell$ factor comes from choosing $B$ and the $\ell(\ell-1)$ factor comes from choosing $B', C'$). To obtain $\sigma|_{B} = \sigma|_{B'}$ and $\sigma|_{C} = \sigma|_{C'}$, we need three subintervals of $k$ vertices to all have the same coloring (since $B$ must have the same coloring as both $B'$ and $C'$), which occurs with probability $(\sum_{i=1}^r \pi_i^3)^k$. The same logic holds for when $B \neq C$ and $B' = C'$. 

The third term comes from the case where both $B \neq C$ and $B' \neq C'$. There are $\ell(\ell-1) \cdot \ell(\ell-1) = \ell^2(\ell-1)^2$ ways to choose such $B, B', C, C'$. Then, $\sigma|_{B} = \sigma|_{B'}$ occurs with probability $q^k$ and $\sigma|_{C} = \sigma|_{C'}$ also occurs with probability $q^k$; since these events are independent, the probability of $\sigma|_{B} = \sigma|_{B'}$ and $\sigma|_{C} = \sigma|_{C'}$ is $q^{2k}$.

Now, we apply Chebyshev's inequality, which yields
\begin{align*}
    \mathbb{P}(Z = 0) &\leq \frac{\E[Z^2]}{\E[Z]^2} - 1 \\
    &\leq \frac{\ell^2q^k + 2\ell^3 (\sum_{i=1}^r \pi_i^3)^k + \ell^4q^{2k}}{\ell^4 q^{2k}} - 1 \\
    &= \frac{1}{\ell^2 q^k} + \frac{2(\sum_{i=1}^r \pi_i^3)^k}{\ell q^{2k}}.
\end{align*}
We will show that both terms go to zero. For the first term, recall that $\ell = n / (6k)$ and observe that $k \leq \frac{2(1-\epsilon)\log n}{\log(1/q)}$ implies $q^k \geq n^{-2(1-\epsilon)}$. Hence, we obtain
\begin{align*}
    \frac{1}{\ell^2q^k} \leq \frac{1}{\frac{n^2}{36k^2} n^{-2(1-\epsilon)}} = \frac{36k^2}{n^{2\epsilon}} \to 0.
\end{align*}
For the second term, we use the monotonicity of $p$-norms \cite{johnston2021advanced} to obtain that $(\sum_i \pi_i^3)^{1/3} \leq q^{1/2}$, which implies that $(\sum_i \pi_i^3) \leq q^{3/2}$. In addition, observe that $k \leq \frac{2(1-\epsilon)\log n}{\log(1/q)}$ implies $q^{k/2} \geq n^{-(1-\epsilon)}$. Thus,
\begin{align*}
    \frac{2(\sum_{i=1}^r \pi_i^3)^k}{\ell q^{2k}} \leq \frac{2q^{3k/2}}{\ell q^{2k}} = \frac{12k}{n q^{k/2}} \leq \frac{12k}{n n^{-(1-\epsilon)}} = \frac{12k}{n^\epsilon} \to 0.
\end{align*}

Therefore, with high probability, there exist $B_1 \in \Gamma_1, B_4 \in \Gamma_4$ such that $\sigma|_{B_1} = \sigma|_{B_4}$. Identical reasoning shows that with high probability, there exist $B_3 \in \Gamma_3, B_6 \in \Gamma_6$ such that $\sigma|_{B_3} = \sigma|_{B_6}$. Hence, \eqref{eq:non-identifiability_DingLiu_1d_condition1} holds with high probability.

Now, we show \eqref{eq:non-identifiability_DingLiu_1d_condition2} using a straightforward calculation and the union bound, following the same reasoning as \cite{ding2023lattice}. For any given $J$, $J'$, we have that
\begin{equation*}
    \bP(\sigma|_{J[1, m]} = \sigma|_{J'[1, m]}) = q^m = q^{n/6},
\end{equation*}
which is exponentially decreasing in $n$. The number of possible choices for $J$ and $J'$ is $\ell^4 = n^4 / (6k)^4$ (coming from $\ell$ choices for each of $B, B', C, C'$), which is polynomially increasing in $n$. By the union bound, \eqref{eq:non-identifiability_DingLiu_1d_condition2} holds with high probability.
\end{proof}

\section{Discussion on the Threshold}
\label{sec:threshold_discussion}
In this section, we propose two possible values for the threshold and discuss why each might be the true threshold. One possible value for the threshold is $k^d = \frac{\log n^d}{H(\pi)}$. There are a couple of reasons for this: First, this value coincides with the threshold in the uniform case, as $H(\pi) = \log r$ in the uniform case. In addition, we have already proven non-identifiability below this value, and despite our attempts, we were unable to extend our proof to larger values of $k^d$. 

However, we do not rule out this possibility. Our proof depends on finding an encoding of the deck such that $H(D) = o(n^d)$ to prove Property 3 in Proposition \ref{prop:2d_negative_fano}, and the threshold is used to show that such an encoding is possible. We think that it may be possible to find a better encoding of the deck such that $H(D) = o(n^d)$ even when $k^d \geq \frac{(1+\epsilon)\log n^d}{H(\pi)}$, which would prove non-identifiability of the lattice for a higher value (and disqualify $k^d = \frac{\log n^d}{H(\pi)}$ as a possibility for the true threshold). 

Another possible value for the threshold is $k^d = \frac{\log n^d}{\log(1/q)}$, where $q = \sum_i \pi_i^2$. This value also coincides with the threshold in the uniform case, as $q = \sum_{i=1}^r (1/r)^2 = 1/r$, so $1/q = r$. Furthermore, in the uniform case, the threshold for $d\geq 2$ is exactly half of the threshold for $d = 1$, which are $k^d = \frac{\log n^d}{\log r}$ and $k^d = \frac{2 \log n^d}{\log r}$ respectively. From \cite{Dyer_Frieze_Suen}, we know that the threshold for the non-uniform case with $d = 1$ is $k^d = \frac{2 \log n^d}{\log(1/q)}$; taking half of this value suggests that the threshold for $d \geq 2$ could be $k^d = \frac{\log n^d}{\log (1/q)}$.

To compare these two possible thresholds, we observe that $H(\pi) \geq \log(1/q)$ with equality holding if and only if $\pi$ is uniform. This can be shown using Jensen's inequality: Let $X$ be the random variable which equals $\pi_i$ with probability $\pi_i$. Using the fact that $f(x) = \log(1/x)$ is convex, we obtain
\begin{equation*}
    H(\pi) = \sum_{i=1}^r \pi_i \log\bigg(\frac{1}{\pi_i}\bigg) = \E\bigg[\log\bigg(\frac{1}{X}\bigg)\bigg] \geq \log\bigg(\frac{1}{\E[X]}\bigg) = \log\bigg(\frac{1}{\sum_{i=1}^r \pi_i^2}\bigg) = \log\bigg(\frac{1}{q}\bigg).
\end{equation*}
Consequently, we see that
\begin{equation*}
    \frac{\log n^d}{H(\pi)} \leq \frac{\log n^d}{\log(1/q)},
\end{equation*}
where equality holds if and only if $\pi$ is uniform. Note that since $H(\pi) = \log(1/q)$ in the uniform case, the difference between these two possible thresholds only appears in the non-uniform case.

There is a qualitative transition in the behavior of the lattice at each of these values for $k^d$. If $k^d \geq \frac{(1+\epsilon)\log n^d}{\log(1/q)}$, then the expected number of times a randomly chosen $k$-card repeats goes to 0.

\begin{lemma}
\label{lem:expected_number_of_repeats}
    If $k^d \geq \frac{(1+\epsilon) \log n^d}{\log(1/q)}$, then the expected number of cards with the same coloring as a randomly chosen $k$-card is $o(1)$.
\end{lemma}
\begin{proof}
    We will consider cards which overlap with the chosen card separately from the cards which do not overlap with the chosen card.

    Observe that the probability that two fixed, non-overlapping cards have the same coloring is $q^{k^d}$, since the probability two vertices have the same color is $q$, and each card has $k^d$ vertices. Thus, the expected number of non-overlapping cards with the same coloring as a randomly chosen card is at most $n^d q^{k^d}$. This is $o(1)$ when $k^d \geq  \frac{(1+\epsilon)\log n^d}{\log(1/q)}$ because $q^{k^d} \leq n^{-(1 + \epsilon)d}$.

    Meanwhile, there are less than $(2k)^d$ cards which overlap the chosen card. Since each card does not overlap with the chosen card on a region of size at least $k^{d-1}$ vertices, the probability that two fixed, overlapping cards have the same coloring is at most $q^{k^{d-1}}$. Thus, the contribution of overlapping cards to the expectation is at most $(2k)^d q^{k^{d-1}} = o(1)$.
\end{proof}

On the other hand, if $\frac{(1+\epsilon) \log n^d}{H(\pi)} \leq k^d \leq \frac{(1-\epsilon)\log n^d}{\log(1/q)}$ (i.e., in between the two proposed values), then a randomly chosen card is unique with high probability, but the expected number of repeats diverges to infinity, as established by the following two lemmas.

\begin{lemma}
    If $k^d \leq \frac{(1-\epsilon)\log n^d}{\log(1/q)}$, then the expected number of cards with the same coloring as a randomly chosen card is $\omega(1)$.
\end{lemma}
\begin{proof}
    The probability that any two fixed cards have the same coloring is at least $q^{k^d}$ (including cards which overlap). Therefore, the expected number of repeats is at least $n^d q^{k^d}$, which tends to infinity when $k^d \leq \frac{(1-\epsilon) \log n^d}{\log(1/q)}$.
\end{proof}

\begin{lemma}
    If $k^d \geq \frac{(1+\epsilon) \log n^d}{H(\pi)}$, then a randomly chosen card is unique with high probability.
\end{lemma}
\begin{proof}
    Fix $\delta = \epsilon /2$. Let $X = (X_1, \ldots, X_r)$, where $X_i$ is the number of vertices of color $i$. We use the law of total probability, conditioning on $X$. Furthermore, we divide the sum into two parts depending on whether the card is $\delta$-typical, in the sense of the event $E_\delta$ defined in \eqref{eq:E_delta}.  Let $A$ be the event that a randomly chosen card is unique. Then,
    \begin{align*}
        \P(A) &= \sum_{x \in E_\delta} \P(A \mid X = x) \P(X = x) + \P(A \cap \{X \not \in E_{\delta}\} ) \\
        &= \sum_{x \in E_\delta} \P(A \mid X = x)\P(X = x) + o(1)
    \end{align*}
    because $E_\delta$ occurs with high probability. Therefore, it is sufficient to show that $\P(A \mid X = x) = 1 - o(1)$ when $x \in E_\delta$ (i.e. that a random typical card is unique with high probability).
 
Recall that the probability that any given card has color count $x$ equals $\prod_{i=1}^r \pi_i^{x_i}$. Under the event $E_\delta$, we have 
\begin{align*}
    \prod_{i=1}^r \pi_i^{x_i} \leq \bigg(\prod_{i=1}^r \pi_i^{(1 - \delta)\pi_i} \bigg)^{k^d}.
\end{align*}
Then, taking a union bound over all $n^d$ cards in the deck, we obtain
\begin{align*}
    \bP(A \mid X = x)  \geq 1 - n^d  \bigg(\prod_{i=1}^r \pi_i^{(1 - \delta)\pi_i} \bigg)^{k^d}.
\end{align*}
We want to show this probability goes to 1. 
Therefore, it is sufficient to show that 
\begin{align*}
    \log(n^d) + k^d \sum_{i=1}^r (1-\delta)\pi_i \log \pi_i \to -\infty.
\end{align*}
Observe that when $k^d \geq \frac{(1+\epsilon) \log n^d}{H(\pi)}$, we have
\begin{align*}
    \log (n^d) + k^d \sum_{i=1}^r (1-\delta)\pi_i \log \pi_i
    &= \log(n^d) - (1-\delta) k^d H(\pi) \\
    &\leq \log(n^d) - (1-\delta)(1+\epsilon) \log (n^d) \\
    &\leq \log(n^d) - (1 + \epsilon/2)\log (n^d) \\
    &= (-\epsilon/2) \log (n^d) \\
    & \to -\infty
\end{align*}
as desired, where the second inequality follows from $\delta = \epsilon / 2$.
\end{proof}

We know that in  the uniform case, if the expected number of times a randomly chosen $k$-card appears goes to zero (as in Lemma \ref{lem:expected_number_of_repeats}) then the lattice is identifiable. We expect the same to hold in the non-uniform case as well, so we believe that the lattice should be identifiable when $k^d \geq \frac{(1+\epsilon) \log n^d}{\log(1/q)}$.

Meanwhile, when $\frac{(1+\epsilon) \log n^d}{H(\pi)} \leq k^d \leq \frac{(1-\epsilon)\log n^d}{\log(1/q)}$, a randomly chosen card is unique with high probability, while the expected number of repeats goes to infinity. Consequently, we have a small number of cards which repeat many times throughout the lattice. This potentially creates a barrier for reconstruction because we might have many regions, centered around these repeated cards, that look similar. This could provide opportunities to swap two cards to create a different lattice but keep the same deck. Therefore, whether the true threshold is $\frac{\log n^d}{H(\pi)}$ or $\frac{\log n^d}{\log(1/q)}$ depends on whether this barrier can be overcome to reconstruct the lattice.

\section{Conclusion and Future Work}
\label{sec:conclusion}
In this paper, we studied the problem of shotgun assembly on the non-uniform labeled lattice model. Our main contribution is establishing a sufficient condition for non-identifiability which matches the known threshold in the uniform case. This work gives rise to some natural open questions.

\begin{enumerate}
    \item What is the true threshold for identifiability? Is it one of the two candidates proposed in Section \ref{sec:threshold_discussion}?
    \item When the lattice is identifiable, what algorithm can be used to reconstruct the lattice?
\end{enumerate}

\bibliographystyle{plain}
\bibliography{refs.bib}

\end{document}